\documentclass[11pt,reqno]{amsart}
\usepackage{graphicx}
\usepackage{verbatim}
\usepackage{textcomp}
\usepackage{amssymb}
\usepackage{cite}
\usepackage{amsmath}
\usepackage{latexsym}
\usepackage{amscd}
\usepackage{amsthm}
\usepackage{mathrsfs}
\usepackage{bm}
\usepackage{url}
\usepackage{hyperref}
\usepackage{bookmark}
\allowdisplaybreaks[3]
\newtheorem{thm}{Theorem}[section]
\newtheorem{corr}[thm]{Corollary}
\newtheorem{lem}[thm]{Lemma}
\newtheorem{prop}[thm]{Proposition}

\theoremstyle{definition}

\newtheorem*{ack}{Acknowledgment}
\theoremstyle{remark}
\newtheorem{rem}{Remark}[section]
\numberwithin{equation}{section}
\begin{document}
\title[Rigidity of Einstein metrics on closed manifolds]
{Rigidity of Einstein metrics as critical points of quadratic curvature functionals on closed manifolds}

\author{Bingqing Ma}
\address{College of Physics and Materials Science, Henan Normal
University, Xinxiang 453007, People's Republic of China}
\address{Department of Mathematics, Henan Normal
University, Xinxiang 453007, People's Republic of China}
\email{bqma@henannu.edu.cn }

\author{Guangyue Huang }
\address{Department of Mathematics, Henan Normal
University, Xinxiang 453007, People's Republic of China}
\email{hgy@henannu.edu.cn}
\email{xxl@henannu.edu.cn}
\email{yuchern@126.com}

\author{Xingxiao Li }

\author{Yu Chen}

\thanks{Research supported by NSFC (Nos. 11401179, 11371018, 11671121).
}

\begin{abstract}
In this paper, we prove some rigidity results for the Einstein metrics as the critical points
of a family of known quadratic curvature functionals on closed manifolds, characterized by some point-wise inequalities. Moreover, we also provide a few rigidity results that involve the Weyl curvature, the trace-less Ricci curvature and the Yamabe invariant, accordingly.
\end{abstract}



\maketitle

\section{Introduction}

In this paper, we always assume that $M^n$ is a closed manifold of dimension $n\geq3$ and $g$ a Riemannian metric on $M^n$ with the Riemannian curvature tensor $R_{ijkl}$, the Ricci tensor $R_{ij}$ and the scalar curvature $R$.
It is well-known that any Einstein metric $g$
must be critical for the Einstein-Hilbert functional
\begin{equation*}
\mathcal {H}=\int_{M}R
\end{equation*}
defined on the space $\mathscr{M}_1(M^n)$ of equivalence classes of
smooth Riemannian metrics of volume one on $M^n$. On the other hand, Catino considered in \cite{Catino2015} the following family of quadratic curvature functionals
\begin{equation}\label{Int-1}
\mathcal{F}_t\,=\,\int_M |R_{ij}|^{2}+t\int_M
R^{2},\quad t\in \mathbb{R}
\end{equation}
which are also defined on $\mathscr{M}_1(M^n)$, and proved some related rigidity results. Furthermore,
it has been observed in \cite{Bess2008} that every
Einstein metric is a critical point of $\mathcal{F}_t$ for all $t\in \mathbb{R}$, see \eqref{2-lemmas-5} below. But the converse of this conclusion is not true in general.

Therefore it is natural to ask that under what conditions a critical metric for the functionals
$\mathcal{F}_t$ must be a Einstein one. In fact, there have been a number of interesting conclusions to this problem, for example, under some suitable curvature conditions (\cite{Lamontagne1994,Catino2015,Tanno1975}), or under some integral conditions (\cite{Gursky2001,HC-2017ProcAMS}). For other development in this direction, we refer the readers to
\cite{HL2004,Gursky2015,Anderson1997,Catino2014ProcAMS}
and the references therein.

Recall that the Yamabe invariant $Y_M([g])$ is defined by
\begin{align}\label{Int-2}
Y_M([g])=&\inf\limits_{\tilde{g}\in[g]
}\frac{\int_M\tilde{R}\,dv_{\tilde{g}}} {(\int_M\,dv_{\tilde{g}})^{\frac{n-2}{n}}}\notag\\
=&\frac{4(n-1)}{n-2}\inf\limits_{u\in W^{1,2}(M^n)}\frac{\int_M|\nabla
u|^2\,dv_g+\frac{n-2}{4(n-1)}\int_M
Ru^2\,dv_g}{(\int_M|u|^{\frac{2n}{n-2}}\,dv_g)^{\frac{n-2}{n}}},
\end{align}
where $[g]$ is the conformal class of the metric $g$. It then follows that
\begin{align}\label{Int-3}
\frac{n-2}{4(n-1)}Y_M([g])& \Big(\int_M|u|^{\frac{2n}{n-2}}\,dv_g\Big)^{\frac{n-2}{n}} \notag\\
\leq&\int_M|\nabla
u|^2\,dv_g+\frac{n-2}{4(n-1)}\int_MRu^2\,dv_g,
\end{align}
for all $u\in W^{1,2}(M^n)$. Moreover, $Y_M([g])$ is positive if and only if there exists a conformal metric in $[g]$ with everywhere positive scalar curvature.

In the present paper, by using some pinching conditions involving the Weyl curvature, the trace-less Ricci curvature and the Yamabe invariant, we aim to prove a number of rigidity theorems for the Einstein metrics considered as the critical points of the functional family $\mathcal{F}_t$ ($t\in R$). For convenience, we shall use $\mathring{{\rm Ric}}$ and $W$ throughout this paper to denote the trace-less Ricci curvature and the Weyl curvature, respectively.

Our main results are stated as follows.

\begin{thm}\label{thm1-1}
Let $(M^n,g)$ be a closed Riemannian manifold of dimension $n\geq 3$ with
positive scalar curvature and $g$ is a critical metric for
the functional family $\mathcal{F}_t$ over $\mathscr{M}_1(M^n)$, where
\begin{equation}\label{1111thm-Int-1}
\begin{cases}t<-\frac{5}{12},\ \ \ \quad \quad {\rm if}\ n=3;\\
t<-\frac{1}{3},\ \quad\quad \quad {\rm if}\ n=4; \\
t\leq-\frac{n}{4(n-1)},\ \quad {\rm if}\ n\geq5.
\end{cases}
\end{equation}
Suppose that
\begin{align}\label{1thm-Int-4}
\Big|W-&\frac{n-4}{\sqrt{2n}(n-2)}\mathring{{\rm Ric}} \mathbin{\bigcirc\mkern-15mu\wedge} g\Big|\notag\\
&<-\sqrt{\frac{2}{(n-1)(n-2)}}\left(\frac{2(n-2)+2n(n-1)t}{n} +1\right)R.
\end{align}
Then $(M^n,g)$ must be of Einstein.
\end{thm}

In particular, when $n=3$, we have $W=0$ automatically. On the other hand, from \eqref{2-lemmas-2} it is seen that
an Einstein manifold $M^3$ with positive scalar curvature must be of constant positive sectional curvature. Moreover, it follows from Lemma \ref{Lemma23} that $|\mathring{{\rm Ric}} \mathbin{\bigcirc\mkern-15mu\wedge} g|=2|\mathring{{\rm Ric}}|$ on $M^3$. Consequently, the following conclusion is immediate by Theorem \ref{thm1-1}.

\begin{corr}\label{corr1-1}
Let $(M^3,g)$ be as in Theorem \ref{thm1-1} with $n=3$. If
\begin{equation}\label{1corr-Int-4}
|\mathring{{\rm Ric}}|<-\frac{5+12t}{\sqrt{6}}R,\quad \text{for}\quad t<-\frac{5}{12},
\end{equation}
then $(M^3,g)$ must be of constant positive sectional curvature.
\end{corr}

Next, for $t=-\frac{1}{2}$, we give rigidity results by using pointwise inequalities.

\begin{thm}\label{thm1-6}
Let $(M^n,g)$ be a closed Riemannian manifold of dimension $n\geq 3$ with
positive scalar curvature and $g$ is a critical metric for
the functional $\mathcal{F}_{-\frac12}$ over $\mathscr{M}_1(M^n)$. Suppose that
\begin{equation}\label{6thm-Int-12}
\Big|W+\frac{\sqrt{2}}{\sqrt{n}(n-2)}\mathring{{\rm Ric}} \mathbin{\bigcirc\mkern-15mu\wedge} g\Big|\leq
\frac{n^2-3n+4}{n\sqrt{2(n-1)(n-2)}} R.
\end{equation}
If there exists a point where the inequality in \eqref{6thm-Int-12} is strict, then
$(M^n,g)$ must be of Einstein.
\end{thm}

To state the next theorem, we first introduce a constant
\begin{equation}\label{2thm-Int-7}
C_n=\begin{cases}\frac1{\sqrt{6}},&\text{if\ }n=4;\\
\frac{3\sqrt{15}-6}{8\sqrt{10}},&\text{if\ }n=5;\\
\frac{2}{n}\left(\frac{2(n-2)}{\sqrt{n(n-1)}} +\frac{n^2-n-4}{\sqrt{n(n-1)(n+1)(n-2)}}\right)^{-1},&\text{if\ }n\geq 6.
\end{cases}
\end{equation}
Then, in terms of the Yamabe invariant, we have the following theorem:

\begin{thm}\label{thm1-2}
Let $(M^n,g)$ be a closed Riemannian manifold of dimension $n\geq 3$ with
positive scalar curvature and $g$ is a critical metric for
the functional $\mathcal{F}_t$ ($t\leq-\frac{1}{2}$) over $\mathscr{M}_1(M^n)$.
Suppose that
\begin{align}\label{2thm-Int-5}
\Big(\int_M\Big|W+\frac{\sqrt{2}}{\sqrt{n}(n-2)}\mathring{{\rm Ric}} \mathbin{\bigcirc\mkern-15mu\wedge} g\Big|^{\frac{n}{2}}\Big)^{\frac{2}{n}}\leq\frac{1}{4}\sqrt{\frac{n-2}{2(n-1)}}Y_M([g]).
\end{align}
Then $(M^n,g)$ must be of Einstein. Furthermore,

(1) if $n=3,4,5$, then \eqref{2thm-Int-5} implies that $(M^n,g)$ is of constant positive sectional curvature;

(2) if $n\geq6$ and \eqref{2thm-Int-5} is replaced with
\begin{align}\label{2thm-Int-6}
\Big(\int_M\Big|W+\frac{\sqrt{2}}{\sqrt{n}(n-2)}\mathring{{\rm Ric}} \mathbin{\bigcirc\mkern-15mu\wedge} g\Big|^{\frac{n}{2}}\Big)^{\frac{2}{n}}<C_nY_M([g]),
\end{align}
then $(M^n,g)$ must be of constant positive sectional curvature.
\end{thm}

\begin{corr}\label{4corr-1}
Let $(M^{4},g)$ be as in Theorem \ref{thm1-2} with $n=4$.
Suppose that
\begin{align}\label{4corr-Int-1}
\int_M\Big(|W|^2+\frac{5}{4}|\mathring{R}_{ij}|^2\Big)
\leq\frac{1}{48}\int_MR^2.
\end{align}
Then, $(M^{4},g)$ must be of constant positive sectional curvature.
\end{corr}

For a four dimensional Riemannian manifold, we also have the following rigidity theorem:

\begin{thm}\label{thm1-3}
Let $(M^{4},g)$ be closed with
positive scalar curvature and $g$ be critical for
$\mathcal{F}_{t}$ over $\mathscr{M}_1(M^{4})$, where
\begin{align}\label{3thm-Int-8}
-\frac{1}{4}\leq t<-\frac{1}{6}.
\end{align}
Suppose that
\begin{align}\label{3thm-Int-9}
\Big(\int_M\Big|W+\frac{1}{2\sqrt{2}}\mathring{{\rm Ric}}
\mathbin{\bigcirc\mkern-15mu\wedge}
g\Big|^{2}\Big)^{\frac{1}{2}}<-\frac{1+6t}{2\sqrt{3}}Y_M([g]).
\end{align}
Then, $(M^{4},g)$ must be of constant positive sectional curvature.
\end{thm}

\begin{corr}\label{3corr-1}
Let $(M^{4},g)$ be as in Theorem \ref{thm1-3}.
If
\begin{align}\label{3corr-Int-1}
\int_M\Big(|W|^2+[1+(1+6t)^2]|\mathring{R}_{ij}|^2\Big)
\leq\frac{(1+6t)^2}{12}\int_MR^2,
\end{align}
then $(M^{4},g)$ must be of constant positive sectional curvature.
\end{corr}

We remark that our rigidity results in both Corollary \ref{4corr-1} and Corollary \ref{3corr-1} can also be described in terms of the Euler-Poincar\'{e} characteristic. In fact, by the well-known Chern-Gauss-Bonnet formula (\cite[Equation 6.31]{Bess2008})
\begin{align}\label{4rem-1}
\int_M\Big(|W^2-2|\mathring{R}_{ij}|^2+\frac{1}{6}R^2\Big)=32\pi^2\chi(M),
\end{align}
with $\chi(M)$ being the Euler-Poincar\'{e} characteristic of $M^4$, the pinching conditions \eqref{4corr-Int-1} and \eqref{3corr-Int-1} are equivalent to
\begin{equation}\label{4rem-2}
\frac{13}{2}\int_M|W|^2+\frac{1}{3}\int_MR^2\leq80\pi^2\chi(M)
\end{equation}
and
\begin{equation}\label{4rem-3}
\frac{3+(1+6t)^2}{2}\int_M|W|^2+\frac{1}{12}\int_MR^2
\leq16[1+(1+6t)^2]\pi^2\chi(M),
\end{equation}
respectively. So we have the following corollary.

\begin{corr}\label{5corr-1}
Let $(M^{4},g)$ be closed with
positive scalar curvature. If either
(1) $g$ is the critical point of
$\mathcal{F}_{t}$ over $\mathscr{M}_1(M^{4})$ with $t\leq-\frac{1}{2}$, and $(M^{4},g)$ satisfies \eqref{4rem-2}, or
(2) $g$ is the critical point of
$\mathcal{F}_{t}$ over $\mathscr{M}_1(M^{4})$ with \eqref{3thm-Int-8}, and $(M^{4},g)$ satisfies \eqref{4rem-3}, then
$(M^{4},g)$ must be of constant positive sectional curvature.
\end{corr}

\begin{rem}\label{1Rem-1}
It should be emphasized that our Corollary \ref{corr1-1} greatly improves
a similar Theorem of Catino in \cite{Catino2015}.

In fact, combining with Remark 1.7 of Catino \cite{Catino2015} and (i) of Theorem 1.1 in \cite{HC-2017ProcAMS}, the mentioned Theorem of Catino can be stated as follows:

\begin{thm}[\cite{Catino2015}, Theorem 1.5] Let $(M^3,g)$ be a Riemannian manifold with
positive scalar curvature and $g$ be a critical metric for
$\mathcal{F}_t$ with $t\in(-\infty, -\frac{1}{2}]\cup [-\varepsilon_0,-\frac{1}{6})$, where $\varepsilon_0\approx 0.3652$. Then $g$ must have constant positive sectional curvature if
\begin{equation}\label{1Rem-Int-3}
|\mathring{{\rm Ric}}|<-\frac{1+6t}{2\sqrt{6}}R.
\end{equation}
\end{thm}

But it is easy to check that
\begin{equation}\label{1Rem-Int-4}
-\frac{1+6t}{2\sqrt{6}}<-\frac{5+12t}{\sqrt{6}}, \quad\text{for\ }t\in(-\infty, -\frac{1}{2}),
\end{equation}
which shows that our pinching condition \eqref{1corr-Int-4} is better than that of Catinos' in this case. Furthermore, the interval $(-\infty,-\frac{5}{12})$ we use for $t$ is clearly larger than the interval $(-\infty, -\frac{1}{2}]\cup [-\varepsilon_0,-\frac{1}{6})$ used by Catinos.
\end{rem}

We should remark that, when $n\geq4$, our Theorem \ref{thm1-1} also generalizes the conclusion (ii) of Theorem 1.1 in \cite{HC-2017ProcAMS}.

\begin{ack}
The second author of this paper thanks Professor Haiping Fu for sending him the recent achievements in this direction.
\end{ack}

\section{Some necessary lemmas}

Recall that the Weyl curvature $W_{ijkl}$ of a Riemannian manifold $(M^n,g)$ with $n\geq3$ is related to the Riemannian curvature $R_{ijkl}$ by
\begin{align}\label{2-lemmas-1}
W_{ijkl}=&R_{ijkl}-\frac{1}{n-2}(R_{ik}g_{jl}-R_{il}g_{jk}
+R_{jl}g_{ik}-R_{jk}g_{il})\notag\\
&+\frac{R}{(n-1)(n-2)}(g_{ik}g_{jl}-g_{il}g_{jk}).
\end{align}
Since the traceless Ricci curvature $\mathring{R}_{ij}=R_{ij}-\frac{R}{n}g_{ij}$,
\eqref{2-lemmas-1} can be written as
\begin{align}\label{2-lemmas-2}
W_{ijkl}=&R_{ijkl}-\frac{1}{n-2}(\mathring{R}_{ik}g_{jl}-\mathring{R}_{il}g_{jk}
+\mathring{R}_{jl}g_{ik}-\mathring{R}_{jk}g_{il})\notag\\
&-\frac{R}{n(n-1)}(g_{ik}g_{jl}-g_{il}g_{jk}).
\end{align}
Furthermore, the Cotton tensor is defined by
\begin{align}\label{2-lemmas-3}
C_{ijk}=&R_{kj,i}-R_{ki,j}-\frac{1}{2(n-1)} (R_{,i}g_{jk}-R_{,j}g_{ik})\notag\\
=&\mathring{R}_{kj,i}-\mathring{R}_{ki,j} +\frac{n-2}{2n(n-1)}(R_{,i}g_{jk}-R_{,j}g_{ik}),
\end{align}
where the indices after a comma denote the covariant derivatives. Then
the divergence of the Weyl curvature tensor is related to the Cotton tensor by
\begin{equation}\label{2-lemmas-4}
-\frac{n-3}{n-2}C_{ijk}=W_{ijkl,l}.
\end{equation}

It has been shown by
Catino in \cite[Proposition 2.1]{Catino2015} that a metric $g$ is critical for $\mathcal{F}_t$ over
$\mathscr{M}_1(M^n)$ if and only if it satisfies the following equations
\begin{align}
\Delta
\mathring{\rm R}_{ij}=&(1+2t)R_{,ij}-\frac{1+2t}{n} (\Delta R)
g_{ij}-2R_{ikjl}\mathring{\rm R}_{kl}\notag\\
&-\frac{2+2nt}{n}R \mathring{\rm R}_{ij} + \frac{2}{n} |\mathring{\rm R}_{ij}|^{2} g_{ij},\label{2-lemmas-5}\\
[n+4(&n-1)t]\Delta R=(n-4) [| R_{ij}|^{2}+t
R^{2}-\lambda],\label{2-lemmas-6}
\end{align}
where $\lambda =\mathcal{F}_t (g)$.

It is easy to see from \eqref{2-lemmas-5} that
\begin{align}\label{2-lemmas-7}
\frac{1}{2}\Delta|\mathring{\rm R}_{ij}|^2=&|\nabla \mathring{\rm R}_{ij}|^2+ \mathring{\rm R}_{ij}\Delta \mathring{\rm R}_{ij}\notag\\
=&|\nabla \mathring{\rm R}_{ij}|^2+(1+2t) \mathring{\rm R}_{ij}R_{,ij}-2R_{ikjl} \mathring{\rm R}_{kl}\mathring{\rm R}_{ij}
-\frac{2+2nt}{n} R|\mathring{\rm R}_{ij}|^2\notag\\
=&|\nabla \mathring{\rm R}_{ij}|^2+(1+2t) \mathring{\rm R}_{ij}R_{,ij}-\frac{2(n-2)+2n(n-1)t}{n(n-1)} R|\mathring{\rm R}_{ij}|^2\notag\\
&+\frac{4}{n-2}\mathring{\rm R}_{ij}\mathring{\rm R}_{jk}\mathring{\rm R}_{ki}-2W_{ikjl} \mathring{\rm R}_{kl}\mathring{\rm R}_{ij}.
\end{align}
Integrating both sides of \eqref{2-lemmas-7} yields
\begin{align}\label{2-lemmas-8}
0=&\int_M|\nabla \mathring{\rm R}_{ij}|^2-\int_M \Big(2W_{ikjl} \mathring{\rm R}_{kl}\mathring{\rm R}_{ij}
-\frac{4}{n-2}\mathring{\rm R}_{ij}\mathring{\rm R}_{jk}\mathring{\rm R}_{ki}\notag\\
&+\frac{2(n-2)+2n(n-1)t}{n(n-1)}R|\mathring{\rm R}_{ij}|^2\Big)
-(1+2t)\int_M\mathring{\rm R}_{ij,j}R_{,i}\notag\\
=&\int_M|\nabla \mathring{\rm R}_{ij}|^2-\int_M \Big(2W_{ikjl} \mathring{\rm R}_{kl}\mathring{\rm R}_{ij}
-\frac{4}{n-2}\mathring{\rm R}_{ij}\mathring{\rm R}_{jk}\mathring{\rm R}_{ki}\notag\\
&+\frac{2(n-2)+2n(n-1)t}{n(n-1)}R|\mathring{\rm R}_{ij}|^2\Big)
-\frac{(n-2)(1+2t)}{2n}\int_M|\nabla R|^2\end{align}
where we have used $\mathring{\rm R}_{ij,j}=\frac{n-2}{2n}R_{,i}$ in the second equality. Thus, we obtain the following result:

\begin{lem}\label{Lemma21}
Let $M^n$ be a closed manifold and $g$ be a critical metric for
$\mathcal{F}_{t}$ on $\mathscr{M}_1(M^n)$. Then
\begin{align}\label{2-lemmas-9}
\int_M|\nabla& \mathring{R}_{ij}|^2=\int_M\Big(
2W_{ijkl}\mathring{R}_{jl}\mathring{R}_{ik}
-\frac{4}{n-2}\mathring{R}_{ij}\mathring{R}_{jk}
\mathring{R}_{ki}\notag\\
&+\frac{2(n-2)+2n(n-1)t}{n(n-1)}R|\mathring{R}_{ij}|^2+\frac{(n-2)(1+2t)}{2n}|\nabla R|^2\Big).
\end{align}

\end{lem}

Next, by combining the relationship
$$R_{ijkl}\mathring{R}_{jl}\mathring{R}_{ik}=W_{ijkl}\mathring{R}_{jl}\mathring{R}_{ik}
-\frac{1}{n(n-1)}R|\mathring{R}_{ij}|^2-\frac{2}{n-2}\mathring{R}_{ij}\mathring{R}_{jk}
\mathring{R}_{ki}$$
with Proposition 4.1 in \cite{Catino2015} (for the case of dimension three, see \cite[Sect. 4]{Gursky2001}), we obtain

\begin{lem}\label{Lemma22}
Let $M^n$ be a closed manifold. Then
\begin{align}\label{2-lemmas-10}
\int_M|\nabla \mathring{R}_{ij}|^2=&\int_M\Big(
W_{ijkl}\mathring{R}_{jl}\mathring{R}_{ik}
-\frac{n}{n-2}\mathring{R}_{ij}\mathring{R}_{jk}
\mathring{R}_{ki}\notag\\
&-\frac{1}{n-1}R|\mathring{R}_{ij}|^2+\frac{(n-2)^2}{4n(n-1)}|\nabla R|^2+\frac{1}{2}|C_{ijk}|^2\Big).
\end{align}

\end{lem}

The next lemma comes from \cite{HM-2016,FuPeng-2017Hokkaido} (for the case of $\lambda=\frac{2}{n-2}$, see \cite{Catino2016}):

\begin{lem}\label{Lemma23}
For every Riemannian manifold $(M^n,g)$ and any $\lambda\in \mathbb{R}$, the following estimate holds
\begin{align}\label{2-lemmas-11}
\Big|-W_{ijkl}\mathring{R}_{jl}&\mathring{R}_{ik}
+\lambda\mathring{R}_{ij}\mathring{R}_{jk}\mathring{R}_{ki}\Big|\notag\\
\leq&
\sqrt{\frac{n-2}{2(n-1)}}\Big(|W|^2+\frac{2(n-2)\lambda^2}{n}|\mathring{R}_{ij}|^2
\Big)^{\frac{1}{2}}|\mathring{R}_{ij}|^2\notag\\
=&\sqrt{\frac{n-2}{2(n-1)}}\Big|W+\frac{\lambda}{\sqrt{2n}}\mathring{{\rm Ric}} \mathbin{\bigcirc\mkern-15mu\wedge} g\Big||\mathring{R}_{ij}|^2.
\end{align}

\end{lem}

The following lemma comes from \cite{FuXiao2017DGA}, or see \cite[Proposition 1.3]{FuXiao2017Monath}:

\begin{lem}\label{Lemma24}
An Einstein manifold $(M^n,g)$ with positive scalar curvature is of constant positive sectional curvature, provided
\begin{equation}\label{2-lemmas-12}
\Big(\int_M|W|^{\frac{n}{2}}\Big)^{\frac{2}{n}}<C_n\,Y_M([g]),
\end{equation}
where $C_n$
is given by \eqref{2thm-Int-7}.

\end{lem}

\section{Proof of the main results}

\subsection{Proof of Theorem \ref{thm1-1}}
Using \eqref{2-lemmas-9} and \eqref{2-lemmas-10}, it is easy to see
\begin{align}\label{thm1-Proof-1}
0=&\int_M\Big[
-W_{ijkl}\mathring{R}_{jl}\mathring{R}_{ik}
-\frac{n-4}{n-2}\mathring{R}_{ij}\mathring{R}_{jk}
\mathring{R}_{ki}\notag\\
&-\frac{1}{n-1}\Big(\frac{2(n-2)+2n(n-1)t}{n}+1\Big)R|\mathring{R}_{ij}|^2\notag\\
&-\frac{n-2}{2n}\Big((1+2t)-\frac{n-2}{2(n-1)}\Big)|\nabla R|^2+\frac{1}{2}|C_{ijk}|^2\Big].
\end{align}
Substituting the estimate \eqref{2-lemmas-11} with $\lambda=-\frac{n-4}{n-2}$ into \eqref{thm1-Proof-1} gives
\begin{align}\label{thm1-Proof-2}
0\geq\int_M\Bigg[&
\left(-\sqrt{\frac{n-2}{2(n-1)}}\Big|W -\frac{n-4}{\sqrt{2n}(n-2)}\mathring{{\rm Ric}} \mathbin{\bigcirc\mkern-15mu\wedge} g\Big|\right.\notag\\
&\left.-\frac{1}{n-1}\Big(\frac{2(n-2)+2n(n-1)t}{n} +1\Big)R\right)|\mathring{R}_{ij}|^2\notag\\
&-\frac{n-2}{2n}\Big((1+2t)-\frac{n-2}{2(n-1)}\Big)|\nabla R|^2+\frac{1}{2}|C_{ijk}|^2\Bigg].
\end{align}
Note that, when $t$ satisfies \eqref{1111thm-Int-1}, $$\frac{2(n-2)+2n(n-1)t}{n}+1<0,\quad (1+2t)-\frac{n-2}{2(n-1)}\leq0.
$$
Then we can use \eqref{1thm-Int-4} to find
\begin{align}\label{thm1-Proof-3}
0\geq&\int_M\Bigg[
\Big[-\sqrt{\frac{n-2}{2(n-1)}}\Big|W-\frac{n-4}{\sqrt{2n}(n-2)}\mathring{{\rm Ric}} \mathbin{\bigcirc\mkern-15mu\wedge} g\Big|\notag\\
&-\frac{1}{n-1}\Big(\frac{2(n-2)+2n(n-1)t}{n}+1\Big)R\Big]|\mathring{R}_{ij}|^2\notag\\
&-\frac{n-2}{2n}\Big((1+2t)-\frac{n-2}{2(n-1)}\Big)|\nabla R|^2+\frac{1}{2}|C_{ijk}|^2\Bigg]
\geq 0,
\end{align}
which shows that $M^n$ is Einstein, completing the proof of Theorem \ref{thm1-1}.

\subsection{Proof of Theorem \ref{thm1-6}}

When $t=-\frac{1}{2}$, we have $(1+2t)R_{,ij}=0$. Thus by \eqref{6thm-Int-12}, the formula \eqref{2-lemmas-7} becomes
\begin{align}\label{thm4-Proof-21}
\frac{1}{2}\Delta|\mathring{\rm R}_{ij}|^2=&|\nabla \mathring{\rm R}_{ij}|^2+\frac{n^2-3n+4}{n(n-1)} R|\mathring{\rm R}_{ij}|^2\notag\\
&+\frac{4}{n-2}\mathring{\rm R}_{ij}\mathring{\rm R}_{jk}\mathring{\rm R}_{ki}-2W_{ikjl} \mathring{\rm R}_{kl}\mathring{\rm R}_{ij}\notag\\
\geq&|\nabla \mathring{\rm R}_{ij}|^2+\Bigg[\frac{n^2-3n+4}{n(n-1)}R\notag\\
&-\sqrt{\frac{2(n-2)}{n-1}}\Big|W+\frac{\sqrt{2}}{\sqrt{n}(n-2)}\mathring{{\rm Ric}} \mathbin{\bigcirc\mkern-15mu\wedge} g\Big|\,\Bigg]|\mathring{R}_{ij}|^2 \geq0
\end{align}
In this case, $|\mathring{\rm R}_{ij}|^2$ is subharmonic on $M^n$. Using the maximum principle, we obtain that $|\mathring{\rm R}_{ij}|$ is constant and $\nabla \mathring{\rm R}_{ij}=0$, implying that the Ricci curvature is parallel, hence the curvature tensor is harmonic and $R$ is constant. In particular, \eqref{thm4-Proof-21} becomes
 \begin{align}\label{thm4-Proof-22}
\Bigg[&\frac{n^2-3n+4}{n\sqrt{2(n-1)(n-2)}} R-\Big|W+\frac{\sqrt{2}}{\sqrt{n}(n-2)}\mathring{{\rm Ric}} \mathbin{\bigcirc\mkern-15mu\wedge} g\Big|\Bigg]|\mathring{R}_{ij}|^2
=0.
\end{align}
If there exists a point $x_0$ such that the inequality \eqref{6thm-Int-12} is strict, then from \eqref{thm4-Proof-22} we have
$|\mathring{R}_{ij}|(x_0)=0$ which with the fact that $|\mathring{\rm R}_{ij}|$ constant shows that $\mathring{R}_{ij}=0$, that is, $M^n$ is Einstein, completing the proof of Theorem \ref{thm1-6}.

\subsection{Proof of Theorem \ref{thm1-2}}

Using the Kato inequality $|\nabla \mathring{R}_{ij}|\geq |\nabla|\mathring{R}_{ij}||$, we have from \eqref{2-lemmas-9}
\begin{align}\label{thm2-Proof-4}
\int_M|\nabla|\mathring{\rm R}_{ij}||^2\leq&\int_M\Big(
2W_{ijkl}\mathring{R}_{jl}\mathring{R}_{ik}
-\frac{4}{n-2}\mathring{R}_{ij}\mathring{R}_{jk}
\mathring{R}_{ki}\notag\\
&+\frac{2(n-2)+2n(n-1)t}{n(n-1)}R|\mathring{R}_{ij}|^2\notag\\ &+\frac{(n-2)(1+2t)}{2n}|\nabla R|^2\Big),
\end{align}
which shows
\begin{align}
\int_M|\nabla|\mathring{\rm R}_{ij}||^2\leq&\frac{(n-2)(1+2t)}{2n}\int_M|\nabla R|^2\notag\\
&+\frac{2(n-2)+2n(n-1)t}{n(n-1)}\int_M R|\mathring{\rm R}_{ij}|^2\notag\\
&+\sqrt{\frac{2(n-2)}{n-1}}\int_M\Big|W +\frac{\sqrt{2}}{\sqrt{n}(n-2)}\mathring{{\rm Ric}} \mathbin{\bigcirc\mkern-15mu\wedge} g\Big||\mathring{R}_{ij}|^2,\label{thm2-Proof-5}
\end{align}
where we have used \eqref{2-lemmas-11}. This together with \eqref{Int-3} gives
\begin{align}\label{thm2-Proof-6}
&\frac{n-2}{4(n-1)}Y_M([g])\Big(\int_M|\mathring{R}_{ij}|^{\frac{2n}{n-2}}\,dv_g\Big)^{\frac{n-2}{n}}\notag\\
\leq&\frac{(n-2)(1+2t)}{2n}\int_M|\nabla R|^2
+\frac{(n+8)(n-2)+8n(n-1)t}{4n(n-1)}\int_M R|\mathring{\rm R}_{ij}|^2\notag\\
&+\sqrt{\frac{2(n-2)}{n-1}}\Big(\int_M\Big|W+\frac{\sqrt{2}}{\sqrt{n}(n-2)}\mathring{{\rm Ric}} \mathbin{\bigcirc\mkern-15mu\wedge} g\Big|^{\frac{n}{2}}\Big)^{\frac{2}{n}}\Big(\int_M|\mathring{R}_{ij}|^{\frac{2n}{n-2}}\,dv_g\Big)^{\frac{n-2}{n}},
\end{align}
where we have used the H\"{o}lder inequality
\begin{align}
&\int_M\Big|W+\frac{\sqrt{2}}{\sqrt{n}(n-2)}\mathring{{\rm Ric}} \mathbin{\bigcirc\mkern-15mu\wedge} g\Big||\mathring{R}_{ij}|^2\notag\\
&\leq\Big(\int_M\Big|W+\frac{\sqrt{2}}{\sqrt{n}(n-2)} \mathring{{\rm Ric}} \mathbin{\bigcirc\mkern-15mu\wedge} g\Big|^{\frac{n}{2}}\Big)^{\frac{2}{n}} \Big(\int_M|\mathring{R}_{ij}|^{\frac{2n}{n-2}}\, dv_g\Big)^{\frac{n-2}{n}}.
\end{align}
In particular, \eqref{thm2-Proof-6} is equivalent to
\begin{align}\label{thm2-Proof-7}
&\Bigg[\frac{n-2}{4(n-1)}Y_M([g])-\sqrt{\frac{2(n-2)}{n-1}}\Big(\int_M\Big|W\notag\\
&+\frac{\sqrt{2}}{\sqrt{n}(n-2)}\mathring{{\rm Ric}} \mathbin{\bigcirc\mkern-15mu\wedge} g\Big|^{\frac{n}{2}}\Big)^{\frac{2}{n}}\Bigg]\Big(\int_M|\mathring{R}_{ij}|^{\frac{2n}{n-2}}\,dv_g\Big)^{\frac{n-2}{n}}\notag\\
\leq&\frac{(n-2)(1+2t)}{2n}\int_M|\nabla R|^2+\frac{(n+8)(n-2)+8n(n-1)t}{4n(n-1)}\int_M R|\mathring{\rm R}_{ij}|^2.
\end{align}
It is easy to check that, for all $n$,
\begin{equation}\label{thm2-Proof-8}
\frac{(n+8)(n-2)}{8n(n-1)}<\frac{1}{2},
\end{equation}
which shows that $1+2t\leq0$ implies $(n+8)(n-2)+8n(n-1)t<0$. Hence, under the assumption \eqref{2thm-Int-5}, $M^n$ must be Einstein.

When $n=4,5$, then \eqref{2thm-Int-5} becomes
\begin{align}\label{thm2-Proof-9}
\Big(\int_M|W|^{\frac{n}{2}}\Big)^{\frac{2}{n}}\leq\frac{1}{4}\sqrt{\frac{n-2}{2(n-1)}}Y_M([g])
\end{align}
since $M^n$ Einstein. By
\begin{align}\label{Pthm2-Proof-10}
\frac{1}{4}\sqrt{\frac{n-2}{2(n-1)}}<C_n,
\end{align}
and with the help of Lemma \ref{Lemma24}, we can derive that $M^n$ is of constant positive sectional curvature, where $C_n$ is given by \eqref{2thm-Int-7}.

When $n\geq6$, we can check that
\begin{align}\label{thm2-Proof-11}
C_n<\frac{1}{4}\sqrt{\frac{n-2}{2(n-1)}}.
\end{align}
Hence, if \eqref{2thm-Int-6} holds, then \eqref{2thm-Int-5} holds. So $M^n$ is of Einstein and hence, by Lemma \ref{Lemma23}, it is also of constant positive sectional curvature, completing the proof of Theorem \ref{thm1-2}.

\subsection{Proof of Theorem \ref{thm1-3}}
In order to prove Theorem \ref{thm1-3}, we shall
need the following proposition.

\begin{prop}\label{3porp-1} Let $M^n$ be a closed manifold with positive scalar curvature and
$g$ be a critical metric for $\mathcal{F}_{t}$ on
$\mathscr{M}_1(M^n)$ with $(1+2t)R_{,ij}=0$. Then, for any
$\alpha>\frac{1}{2}$, we have
\begin{align}\label{3Prop-1}
0\geq&\Bigg[\Big(2-\frac{1}{\alpha}\Big) \frac{n-2}{4(n-1)}Y_M([g])\notag\\
&-\alpha\sqrt{\frac{2(n-2)}{n-1}} \Big(\int_M\Big|W+\frac{\sqrt{2}} {\sqrt{n}(n-2)}\mathring{{\rm
Ric}} \mathbin{\bigcirc\mkern-15mu\wedge}
g\Big|^{\frac{n}{2}}\Big)^{\frac{2}{n}}
\Bigg]\notag\\
&\times\Big(\int_M|\mathring{R}_{ij}|^{\frac{2n\alpha}{n-2}} \Big)^{\frac{n-2}{n}}
-\Bigg[\frac{2\alpha[(n-2)+n(n-1)t]}{n(n-1)}\notag\\
&+\Big(2-\frac{1}{\alpha}\Big)\frac{n-2}{4(n-1)}\Bigg] \int_MR|\mathring{\rm
R}_{ij}|^{2\alpha}.
\end{align}
\end{prop}

\proof  Since $(1+2t)R_{,ij}=0$, \eqref{2-lemmas-7} becomes
\begin{align}\label{3th-Proof1}
|\mathring{\rm R}_{ij}|\Delta|\mathring{\rm R}_{ij}|\geq& -\frac{2(n-2)+2n(n-1)t}{n(n-1)} R|\mathring{\rm R}_{ij}|^2\notag\\
&+\frac{4}{n-2} \mathring{\rm R}_{ij} \mathring{\rm R}_{jk} \mathring{\rm R}_{ki}-2W_{ikjl}  \mathring{\rm R}_{kl} \mathring{\rm R}_{ij}\notag\\
\geq&-\sqrt{\frac{2(n-2)}{n-1}}\Big|W+\frac{\sqrt{2}}{\sqrt{n}(n-2)}\mathring{{\rm Ric}} \mathbin{\bigcirc\mkern-15mu\wedge} g\Big||\mathring{R}_{ij}|^2\notag\\
&-\frac{2(n-2)+2n(n-1)t}{n(n-1)} R|\mathring{\rm R}_{ij}|^2.
\end{align}
Let $u=|\mathring{\rm R}_{ij}|$. Then for any $\alpha>0$, we have
\begin{align}\label{3th-Proof2}
u^{\alpha}\Delta u^{\alpha}=&u^{\alpha}[\alpha(\alpha-1)u^{\alpha-2}|\nabla u|^2+\alpha u^{\alpha-1}\Delta u]\notag\\
=&\Big(1-\frac{1}{\alpha}\Big)|\nabla u^{\alpha}|^2+\alpha u^{2\alpha-2}u\Delta u\notag\\
\geq&\Big(1-\frac{1}{\alpha}\Big)|\nabla u^{\alpha}|^2-\frac{2\alpha[(n-2)+n(n-1)t]}{n(n-1)} Ru^{2\alpha}\notag\\
&-\alpha\sqrt{\frac{2(n-2)}{n-1}}\Big|W+\frac{\sqrt{2}}{\sqrt{n}(n-2)}\mathring{{\rm
Ric}} \mathbin{\bigcirc\mkern-15mu\wedge} g\Big|u^{2\alpha},
\end{align}
and hence
\begin{align}\label{3th-Proof3}
0\geq&\Big(2-\frac{1}{\alpha}\Big)\int_M|\nabla u^{\alpha}|^2-\frac{2\alpha[(n-2)+n(n-1)t]}{n(n-1)} \int_MRu^{2\alpha}\notag\\
&-\alpha\sqrt{\frac{2(n-2)}{n-1}}\int_M\Big|W+\frac{\sqrt{2}}{\sqrt{n}(n-2)}\mathring{{\rm
Ric}} \mathbin{\bigcirc\mkern-15mu\wedge} g\Big|u^{2\alpha}.
\end{align}
Therefore, by virtue of \eqref{Int-3}, we have for $2-\frac{1}{\alpha}>0$

\begin{align}\label{3th-Proof4}
0\geq&\Bigg[\Big(2-\frac{1}{\alpha}\Big) \frac{n-2}{4(n-1)}Y_M([g])\notag\\
&-\alpha\sqrt{\frac{2(n-2)}{n-1}} \Big(\int_M\Big|W+\frac{\sqrt{2}} {\sqrt{n}(n-2)}\mathring{{\rm
Ric}} \mathbin{\bigcirc\mkern-15mu\wedge}
g\Big|^{\frac{n}{2}}\Big)^{\frac{2}{n}}
\Bigg]\notag\\
&\times\Big(\int_Mu^{\frac{2n\alpha}{n-2}}\Big) ^{\frac{n-2}{n}}
-\Bigg[\frac{2\alpha[(n-2)+n(n-1)t]}{n(n-1)}\notag\\
&+\Big(2-\frac{1}{\alpha}\Big)\frac{n-2} {4(n-1)}\Bigg]\int_MRu^{2\alpha},
\end{align}
where we have used the H\"{o}lder inequality
\begin{align}
\int_M\Big|W+&\frac{\sqrt{2}}{\sqrt{n}(n-2)}\mathring{{\rm Ric}} \mathbin{\bigcirc\mkern-15mu\wedge} g\Big|u^{2\alpha}\notag\\
\leq&\Big(\int_M\Big|W+\frac{\sqrt{2}} {\sqrt{n}(n-2)}\mathring{{\rm
Ric}} \mathbin{\bigcirc\mkern-15mu\wedge}
g\Big|^{\frac{n}{2}}\Big)^{\frac{2}{n}}\Big(\int_M
u^{\frac{2n\alpha}{n-2}}\Big)^{\frac{n-2}{n}}.
\end{align}
So the proof of Proposition \ref{3porp-1} is completed. \endproof

Now, we are in a position to prove Theorem \ref{thm1-3}.

For $n=4$, using \eqref{2-lemmas-6}, we obtain that the scalar
curvature $R$ is harmonic and hence
$R$ is constant. Therefore, \eqref{3Prop-1} becomes
\begin{align}\label{thm3-Proof-16}
0\geq&\Bigg[\frac{1}{6}\Big(2-\frac{1}{\alpha}\Big) Y_M([g])-\frac{2\alpha}{\sqrt{3}}\Big(\int_M\Big|W
+\frac{1}{2\sqrt{2}}\mathring{{\rm Ric}} \mathbin{\bigcirc\mkern-15mu\wedge} g\Big|^{2}\Big)^{\frac{1}{2}}
\Bigg]\notag\\
&\times\Big(\int_M|\mathring{R}_{ij}|^{4\alpha}\Big) ^{\frac{1}{2}}
-\Bigg[\frac{\alpha(1+6t)}{3}+\frac{1}{6} \Big(2-\frac{1}{\alpha}\Big)\Bigg]
\int_MR|\mathring{\rm R}_{ij}|^{2\alpha}.
\end{align}
Since $t$ satisfies \eqref{3thm-Int-8}, we may take
$$\alpha=\frac{1-\sqrt{1+2(1+6t)}}{-2(1+6t)}$$
such that $2-\frac{1}{\alpha}>0$. In this case, we have $$\frac{\alpha(1+6t)}{3}+\frac{1}{6} \Big(2-\frac{1}{\alpha}\Big)=0.$$
Hence \eqref{thm3-Proof-16} becomes
\begin{align}\label{thm3-Proof-17}
0\geq&\Bigg[\frac{1}{6}\Big(2-\frac{1} {\alpha}\Big)Y_M([g])-\frac{2\alpha} {\sqrt{3}}\Big(\int_M\Big|W
+\frac{1}{2\sqrt{2}}\mathring{{\rm Ric}} \mathbin{\bigcirc\mkern-15mu\wedge} g\Big|^{2}\Big)^{\frac{1}{2}}
\Bigg]\notag\\
&\times\Big(\int_M|\mathring{R}_{ij}|^{4\alpha}\Big) ^{\frac{1}{2}}.
\end{align}
Therefore, under the assumption \eqref{3thm-Int-9}, we obtain that
$M^4$ is Einstein. Moreover, in this case, \eqref{3thm-Int-9}
becomes
\begin{align}\label{thm3-Proof-19}
\Big(\int_M|W|^{2}\Big)^{\frac{1}{2}}<-\frac{1+6t}{2\sqrt{3}}Y_M([g]).
\end{align}
We can check that
\begin{align}\label{thm3-Proof-20}
-\frac{1+6t}{2\sqrt{3}}<\frac{1}{\sqrt{6}},
\end{align}
which, combined with Lemma \ref{Lemma24}, shows that $M^4$ is of
constant positive sectional curvature, completing
the proof of Theorem \ref{thm1-3}.

\section{Proof of Corollaries \ref{4corr-1} and \ref{3corr-1}}

In this last section, we provide the detail in proving Corollaries \ref{4corr-1} and \ref{3corr-1}. For this, the following lemma by Catino (Lemma 4.1, \cite{Catino2016}) is needed:

\begin{lem}\label{Appendix-Lemma-1}
Let $(M^4,g)$ be a closed  manifold. Then
\begin{align}\label{Appendix-1}
Y_M^2([g])\geq\int_M(R^2-12|\mathring{R}_{ij}|^2),
\end{align}
and the inequality is strict unless $(M^4,g)$  is conformally Einstein.

\end{lem}

When $n=4$, the pinching condition \eqref{2thm-Int-5} can be written as
\begin{align}\label{Appendix-2}
&\int_M(|W|^2+|\mathring{R}_{ij}|^2)
<\frac{1}{48}Y_M^2([g]).
\end{align}
It holds by \eqref{Appendix-1} that
\begin{align}\label{Appendix-3}
&\int_M(|W|^2+|\mathring{R}_{ij}|^2)-\frac{1}{48}Y_M^2([g])\notag\\
<&\int_M\Big(|W|^2+\frac{5}{4}|\mathring{R}_{ij}|^2-\frac{1}{48}R^2\Big),
\end{align}
from which Corollary \ref{4corr-1} follows immediately.

On the other hand, \eqref{3thm-Int-9} can be written as
\begin{align}\label{Appendix-4}
\int_M(|W|^2+|\mathring{R}_{ij}|^2)
<\frac{(1+6t)^2}{12}Y_M^2([g])
\end{align}
which, combined with \eqref{Appendix-1}, gives
\begin{align}\label{Appendix-5}
&\int_M(|W|^2+|\mathring{R}_{ij}|^2)-\frac{(1+6t)^2}{12}Y_M^2([g])\notag\\
<&
\int_M\Big(|W|^2+[1+(1+6t)^2]|\mathring{R}_{ij}|^2-\frac{(1+6t)^2}{12}R^2\Big),
\end{align}
completing the proof of Corollary \ref{3corr-1}.

\bibliographystyle{Plain}

\begin{thebibliography}{10}

\bibitem{Anderson1997}
M.  Anderson, Extrema of curvature functionals on the space of
metrics on 3-manifolds, Calc. Var. Partial Differential Equations, \textbf{5}
(1997), 199-269.

\bibitem{Bess2008} A. Besse, Einstein Manifolds, Springer-Verlag,
Berlin, 2008.

\bibitem{Catino2014ProcAMS}
G. Catino, Critical metrics of the $L^2$-norm of the scalar curvature, Proc. Amer. Math. Soc. \textbf{142} (2014), 3981--3986.

\bibitem{Catino2015} G. Catino, Some rigidity results on critical metrics for quadratic
functionals, Calc. Var. Partial Differential Equations, \textbf{54} (2015),
 2921--2937.

\bibitem{Catino2016} G. Catino,
Integral pinched shrinking Ricci solitons,
Adv. Math. \textbf{303} (2016), 279--294.


\bibitem{FuXiao2017DGA}
H. Fu, L. Xiao,
Einstein manifolds with finite $L^p$-norm of the Weyl curvature, Differential Geom. Appl. \textbf{53} (2017), 293--305.

\bibitem{FuXiao2017Monath}
H. Fu, L. Xiao, Rigidity theorem for integral pinched shrinking Ricci
solitons, Monatsh. Math. \textbf{183} (2017), 487--494.


\bibitem{FuPeng-2017Hokkaido}
H. Fu, J. Peng,
Rigidity theorems for compact Bach-flat manifolds with positive
constant scalar curvature, arXiv:1707.07236, to appear in {\it Hokkaido Math. J.}

\bibitem{Gursky2001}
M. Gursky, J. Viaclovsky, A new variational characterization of
three-dimensional space forms, Invent. Math. \textbf{145} (2001), 251--278.

\bibitem{Gursky2015}
M. Gursky, J. Viaclovsky, Rigidity and stability of Einstein metrics
for quadratic curvature functionals, J. Reine Angew. Math. \textbf{700}
(2015), 37--91.

\bibitem{HL2004}  Z. Hu,  H. Li, A new variational characterization of
$n$-dimensional space forms, Trans. Amer. Math. Soc. \textbf{356} (2004),
3005--3023.

\bibitem{HM-2016Colloq}
G. Huang, B. Ma, Riemannian manifolds with harmonic curvature, Colloq. Math. \textbf{145} (2016), 251--257.

\bibitem{HC-2017ProcAMS}
G. Huang, L. Chen, Some characterizations on critical metrics for quadratic curvature
functions, Proc. Amer. Math. Soc. \textbf{146} (2018), 385--395.

\bibitem{HM-2016}
G. Huang, Rigidity of Riemannian manifolds with positive scalar curvature, arXiv:1707.00902.


\bibitem{Lamontagne1994}
F. Lamontagne, Une remarque sur la norme $L^2$ du tenseur de
courbure, C. R. Acad. Sci. Paris S\'{e}r. I Math. \textbf{319} (1994),
237--240.


\bibitem{Tanno1975}
S. Tanno, Deformations of Riemannian metrics on 3-dimensional
manifolds, T\^{o}hoku Math. J.  \textbf{27} (1975), 437--444.



\end{thebibliography}

\end{document}